\documentclass[12pt]{article}
\usepackage{amsmath}
\usepackage{amsthm}
\usepackage{amssymb}
\usepackage{color}
\usepackage{indentfirst}

\begin{document}

\newtheorem{theorem}{Theorem}
\newtheorem{lemma}{Lemma}
\newtheorem{corollary}{Corollary}
\newtheorem{proposition}{Proposition}
\newtheorem{definition}{Definition}

\theoremstyle{remark}

\newcommand{\wind}{\mathop{\mathrm{wind}}}

\title {Quasisymmetrically minimal homogeneous perfect sets\thanks{This
work is supported by NNSF No.11071059}}
\author{Yingqing Xiao \thanks{College of Mathematics and Economics, Hunan University,
Changsha 410082, China}}

\maketitle

\renewcommand{\thefootnote}{}
\footnote[1]{{\it E-mail address:} ouxyq@yahoo.cn (Yingqing Xiao)}

\textbf{Abstract:}In \cite{ZW}, the notion of homogenous perfect set as a generalization of Cantor type sets is introduced. Their Hausdorff, lower box-counting, upper box-counting and packing dimensions are studied in \cite{ZW} and \cite{WW}. In this paper, we show that the homogenous perfect set be minimal for $1$-dimensional quasisymmetric maps, which generalize the conclusion in \cite{MS} about the uniform Cantor set to the homogenous perfect set.


\smallskip

\textbf{2000 mathematics classification:} Primary 30C62; Secondary 28A78.

\section{Introduction}


Given $M\geq 1$, a homeomorphism $f:\mathbb{R}\rightarrow \mathbb{R}$ is said to be $M-$quasisymmetric
if and only if
$$
M^{-1}\leq \frac{|f(I)|}{|f(J)|}\leq M
$$
for all pairs of adjacent intervals $I,J$ of equal length, here and in sequel $|\cdot|$ stands for the 1-dimensional Lebesgue measure. A
map is quasisymmetric if it is $M-$quasisymmetric for some $M\geq1$. More generally a homeomorphism between metric
spaces $(X,d_X)$ and $(Y,d_Y)$. If there is a homeomorphism
$\eta:[0,+\infty)\rightarrow [0,+\infty)$ such that
\begin{equation}
\frac{d_X(a,x)}{d_X(b,x)}\leq t\Rightarrow
\frac{d_Y(f(a),f(x))}{d_Y(f(b),f(x))}\leq
\eta(t)
\end{equation}
for all triples $a, b, x $ of distinct points in $X$ and $t\in [0,+\infty)$, then we call $f$ is a quasisymmetric map. When $X =Y=\mathbb{R}^n$, we also say that $f$ is an $n$-dimensional quasisymmetric map.


Let $QS(X)$ denote the collection of all quasisymmetric maps defined on $X$.
Conformal dimension of a metric space, a concept introduced by Pansu
in \cite{PP}, is the infimal Hausdorff dimension of quasisymmetric images of $X$,
$$
\mathcal{C}\dim X =\inf_{f\in QS(X)}\dim_Hf(X).
$$
We say $X$ is minimal for conformal dimension or just minimal if $\mathcal{C}\dim X =
\dim_H X$. Euclidean spaces with standard metric are the simplest examples
of minimal spaces. Basic analytic definitions and results about the conformal dimension and the quasisymmetric map are contained in \cite{JH}.


Now, we introduce the notion of the homogeneous perfect set. The general references on the homogeneous perfect set are \cite{ZW,WW}. In these paper, the authors obtained the  Hausdorff, lower box-counting, upper box-counting and packing dimensions of  the homogeneous perfect set.

{\bf Homogeneous perfect sets.} Let $J_0=[0,1]\subset R$ be the fixed closed interval which
we call the initial interval. Let $\{n_k\}^{\infty}_{k =1}$ be a
sequence of positive integers and $\{c_k\}$ a
sequence of positive real numbers such that for any $k\geq 1,n_k\geq 2$ and $0<c_{k}<1$. For any $k\geq 1$, let
$D_k=\{(i_1,i_2,\cdot\cdot\cdot,i_k):1\leq i_j\leq n_j, 1\leq j\leq
k\} , D=\bigcup_{k\geq 0}D_k$, where $D_0=\{0\}$. We assume if
$\sigma=(\sigma_1, \sigma_2,\cdot\cdot\cdot,\sigma_k)\in D_k, 1\leq
j \leq n_{k+1}$, then $\sigma*j= (\sigma_1,
\sigma_2,\cdot\cdot\cdot,\sigma_k,j)\in D_{k+1}$.

Suppose that $J_0$ is the initial interval and $\mathcal{J}=\{J_{\sigma}:\sigma\in D\}$ is a collection of closed subintervals of $J_0$.  We say that the collection $\mathcal{J}$ fulfills the homogenous perfect structure provided:

1. For any $k\geq 0, \sigma\in D_k, J_{\sigma*1}, J_{\sigma*2},\cdot\cdot\cdot, J_{\sigma*n_{k+1}}$ are subintervals of $J_{\sigma}$. Furthermore, $\max\{x:x\in J_{\sigma*i}\}\leq \min\{x:x\in J_{\sigma*(i+1)}\}, 1\leq i\leq n_{k+1}-1$, that is the interval $J_{\sigma*i}$ is located at the left of $J_{\sigma*(i+1)}$ and the interiors of the intervals $J_{\sigma*i}$ and $J_{\sigma*(i+1)}$ are disjoint.

2. For any $k\geq 1, \sigma\in D_{k-1}, 1\leq j\leq n_k$, we have
$$
\frac{|J_{\sigma*i}|}{|J_{\sigma}|}=c_k.
$$

3. There exists a sequence of nonnegative real numbers $\{\eta_{k,j},k\geq 1, 0\leq j\leq n_k\}$ such that for any $k\geq 0, \sigma\in D_k$, we have $\min(J_{\sigma*1})-\min(J_{\sigma})=\eta_{k+1,0},
\max(J_{\sigma})-\max(J_{\sigma*n_{k+1}})=\eta_{k+1,n_{k+1}}$, and $\min(J_{\sigma*(i+1)})-\max(J_{\sigma*i})=\eta_{k+1,i}(1\leq i\leq n_{k+1}-1)$.

Suppose that the collection of intervals $\mathcal{J}=\{J_{\sigma}:\sigma\in D\}$ satisfies the homogeneous perfect structure.

Let $$E_k=\bigcup_{\sigma\in D_k} J_{\sigma}$$ for every $k\geq 1$.
The set
$$
E:= E(J_0,\{n_k\},\{c_k\},\{\eta_{k,j}\})=\bigcap_{k\geq 1}\bigcup_{\sigma\in D_k}J_{\sigma}=\bigcap_{k\geq 0}E_k
$$
is called a homogeneous perfect set and the intervals $J_{\sigma},\sigma\in D_k$, the fundamental intervals of order $k$.

 For any $k\geq 1$, if $\eta_{k,0}=\eta_{k,n_k}=0$ and $\eta_{k,l}=e_k|J_{\sigma}|$ for all $1\leq l\leq n_{k}-1,\sigma\in D_{k-1}$. Then $E$ is called a uniform Cantor set. This case has been considered by M.D. Hu and S.Y.Wen in \cite{MS}. They obtained
\begin{theorem}[\cite{MS}]\label{MSThm}
Let $E$ be a uniform Cantor set. If the sequence $\{n_k\}$ is bounded and if
$\dim_HE=1$. Then $\dim_H
f(E)=1$ for all $1$-dimensional quasisymmetric maps $f$.
\end{theorem}

In this paper, we generalize Theorem \ref{MSThm} to the homogeneous perfect set and show how the techniques of \cite{MS} can be applied to the homogeneous perfect set and obtain the following theorem.
\begin{theorem}\label{MTHM}
Let $E:= E(J_0,\{n_k\},\{c_k\},\{\eta_{k,j}\})$ be a homogeneous perfect set. If the
sequence $\{n_k\}$ is bounded and if $\dim_HE=1$, then
$\dim_Hf(E)=1$ for all $1$-dimensional quasisymmetric map $f$.
\end{theorem}

This paper is organized as following. In section 2 we introduce the
basic general definitions and results in fractal geometry. The proof of Theorem \ref{MTHM} appears in section 3.


\section{Preliminary}
In order to obtain our result, we need the following lemma from \cite{JMW}, the lemma can also be found in \cite{HH} or \cite{MS}.
\begin{lemma}[\cite{JMW}]\label{mainl2}
Let $f$ be an $M$-quasisymmetric map. Then
\begin{equation}\label{numpq1}
(1+M)^{-2}(\frac{|J|}{|I|})^q\leq\frac{|f(J)|}{|f(I)|}\leq
4(\frac{|J|}{|I|})^p
\end{equation}
for all pairs $J,I$ of intervals with $J\subset I$, where
\begin{equation}\label{numpq}
0<p=\log_2(1+M^{-1})\leq 1 \leq q=\log_2(1+M).
\end{equation}
\end{lemma}

{\bf Hausdorff dimension.}
In this subsection, we recall the definition of Hausdorff dimension. For more details we refer to \cite{KF,ZYW}.

Let $K\subset \mathbb{R}^d$. For any $s\geq 0$, the $s-$dimensional Hausdorff measure of $K$ is given in the usual way by
$$
\mathbf{H}^s(K)=\liminf_{\delta\rightarrow 0}\{\sum_i|U_i|^s: K\subset \bigcup_iU_i, 0<|U_i|<\delta\}.
$$
This leads to the definition of the
Hausdorff dimension of $K$:
$$
\dim_HK=\inf\{s: \mathbf{H}^s(K)<\infty\}=\sup\{s: \mathbf{H}^s(K)>0\}.
$$

The Hausdorff dimension of  the homogeneous perfect set $E$, which depends on $\{n_k\},\{c_k\}$ and $\{\eta_{k,j}\}$ have been obtained in \cite{ZW} as follows
\begin{theorem}[\cite{ZW}]\label{Hdim}
Let $E = E(J_0,\{n_k\},\{c_k\},\{\eta_{k,j}\})$ be a homogeneous perfect set. Suppose
$n_k\leq D$ for all $k$, where $D$ is a constant, then
\begin{equation}
\dim_HE=\liminf_{k\rightarrow \infty}\frac{\log(n_{1}n_{2}\cdot\cdot\cdot n_{k})}{-\log(\sum_{l=1}^{n_{k+1}-1}\eta_{k+1,l}+n_{k+1}c_1c_2\cdot\cdot\cdot c_{k+1})}.
\end{equation}
\end{theorem}

Denote by $N_k$ the number of component intervals of $E_k$ and by $\delta_k$ their common length. Let $e_{k,l}=\eta_{k,l}/\delta_{k-1}\geq \eta_{k,l}$ for all $k\geq 1$ and $0\leq l\leq n_{k}$.
From the definition we obtain
$$
n_kc_k\leq 1,\quad
N_k=n_kn_{k-1}\cdot\cdot\cdot n_1 \quad \mathrm{and} \quad \delta_k=c_kc_{k-1}\cdot\cdot\cdot c_1
$$
for all $k\geq1$. So we have the total length of $E_k$ is
$$
N_k\delta_k=\prod_{i=1}^kn_ic_i,
$$
and
\begin{equation}\label{delata}
\delta_k=\Sigma_{l=0}^{n_{k+1}}\eta_{k+1,l}+n_{k+1}\delta_{k+1}
=\Sigma_{l=0}^{n_{k+1}}e_{k+1,l}\delta_k+n_{k+1}\delta_{k+1}.
\end{equation}

\begin{lemma}\label{Mlema}
Let $E = E(J_0,\{n_k\},\{c_k\},\{\eta_{k,j}\})$ be a homogeneous perfect set. Suppose the sequence $\{n_k\}$ is bounded and $\mathrm{dim}_{H}E=1$ then:

$(1)$ $\lim_{k\rightarrow \infty} (N_k\delta_k)^{1/k}=1$.

$(2)$ $\lim_{k\rightarrow \infty} \frac{1}{k}\sum_{i=1}^ke_i^p=0$ for
any $0<p\leq 1$, where $e_i=\max_{0\leq l\leq n_i}{e_{i,l}}$.

$(3)$ $\lim_{k\rightarrow \infty}\sharp\{i: 0\leq i \leq k, e_i\geq
\epsilon\}/k=0$ for any $\epsilon\in (0,1)$,
where $\sharp$ denotes the cardinality.

\end{lemma}
\proof $(1)$  Since
$$N_k(\delta_k-\eta_{k,0}-\eta_{k,n_{k+1}})\leq  N_k\delta_k\leq 1,$$
Thus, we have
$$\frac{\log N_k}{-\log(\delta_k-\eta_{k,0}-\eta_{k,n_{k+1}})}\leq \frac{\log N_k}{-\log\delta_k}\leq 1.$$
As $\mathrm{dim}_{H}E=1$, we get from Theorem \ref{Hdim}
\begin{equation}
\begin{split}
1=\dim_HE &\;=\liminf_{k\rightarrow \infty}\frac{\log N_k}{-\log(\delta_k-\eta_{k,0}-\eta_{k,n_{k+1}})}\\
&\;\leq\lim_{k\rightarrow \infty}\frac{\log N_k}{-\log\delta_k}\leq 1.
\end{split}
\end{equation}
Thus we obtain
$$
\lim_{k\rightarrow \infty}\frac{\log N_k}{-\log\delta_k}=\lim_{k\rightarrow \infty}\frac{\log N_k}{\log N_k-\log N_k\delta_k}=1,
$$
and
$$
\lim_{k\rightarrow \infty}\frac{\log N_k\delta_k}{\log N_k}=0.
$$
Let $N=1+\sup_k n_k<\infty$. We obtain $N_k\leq N^k$,
so
$$
\lim_{k\rightarrow \infty}\frac{\log N_k\delta_k}{k\log N}=0,
$$
that gives the the conclusion $(1)$ of the lemma.

$(2)$ Since
$$(N_k\delta_k)^{1/k}=(\prod_{i=1}^kn_ic_i)^{1/k}\leq \frac{1}{k}\sum_{i=1}^kn_ic_i\leq 1.
$$
Thus, we have
\begin{equation}\label{eq11}
\lim_{k\rightarrow \infty}\frac{1}{k}\sum_{i=1}^kn_ic_i=1.
\end{equation}

From the equation (\ref{delata}), we have
\begin{equation}
\delta_k
=\Sigma_{l=0}^{n_{k+1}}e_{k+1,l}\delta_k+n_{k+1}c_{k+1}\delta_{k}.
\end{equation}
Thus
$$
e_{k+1}\leq  1-n_{k+1}c_{k+1},
$$
so
$$
\frac{1}{k}\sum_{i}^ke_i\leq\frac{1}{k}\sum^k_{i}(1-n_ic_{i}).
$$
Since the equation (\ref{eq11}), we obtain
$$
\lim_{i}\frac{1}{k}\sum_{i}^ke_i=0,
$$
which together with Jensen's
inequality yields
$$
\lim_{k\rightarrow \infty} \frac{1}{k}\sum_{i=1}^ke_i^p\leq
\lim_{k\rightarrow \infty} (\frac{1}{k}\sum_{i=1}^ke_i)^p=0
$$
for any $0<p\leq 1$. This proves the conclusion $(2)$.

$(3)$ Fixed $\epsilon\in (0,1)$, we obtain from the conclusion $(2)$

$$
\frac{1}{k}\sharp\{i: 0\leq i \leq k, e_i\geq \epsilon
\}=\frac{1}{k}\sum_{i:1\leq i\leq k,e_i\geq \epsilon}1\leq
\frac{1}{k\epsilon}\sum_{i=1}^ke_i\rightarrow 0
$$
as $k$ tends to $\infty$. This proves the conclusion $(3)$.

\section{The proof of Theorem \ref{MTHM}}
In order to obtain our result, we need the following mass distribution principle to estate the lower bound.
\begin{lemma}[\cite{KF}]\label{MSSLemma}Let $\mu$ be a mass distribution supported on $E$. Suppose that for some $t$ there are numbers $c>0$ and $\eta>0$ such that for all sets $U$ with $|U|\leq \eta$ we have $\mu(U)\leq c|U|^t$. Then $\dim_HE\geq t.$
\end{lemma}

{\bf The proof of Theorem \ref{MTHM}:}
Let $E=\bigcap_{k=0}^{\infty}E_k$ be a
homogeneous perfect set satisfying the conditions of Theorem \ref{MTHM}. Let $f:
\mathbb{R}\rightarrow \mathbb{R}$ be an $M$-quasisymmetric map and
$q$ is the number defined as in $(\ref{numpq})$.
Without loss of generality assume that $f([0,1])
=[0,1]$. Then $f(E)=\bigcap_{k=1}^{\infty}f(E_k)$.
The images of component intervals of $E_k$ are component intervals
of $f(E_k)$.

We define a mass distribution $\mu$ on $f(E)$ as follows: Let
$\mu([0, 1]) = 1$. For every $k\geq 1$ and for every component
interval $J$ of $f (E_{k-1})$, let $J_{k1},J_{k2}, \cdot\cdot\cdot,
J_{kn_{k}} $denote the $n_k $ component intervals of $f(E_k)$
lying in $J$. Define
$$
\mu(J_{ki})=\frac{|J_{ki}|^{d}}{||J||_{d}}\mu(J),\quad
i=1,2,\cdot\cdot\cdot,n_k,
$$
where
$$
||J||_d=\sum_{i=1}^{n_k}|J_{ki}|^d
$$
and
\begin{equation} d\in\left\{
\begin{array}{ll}
(0,1)& $ when $ q=1,\\
(1/q,1)& $ when $q>1.
\end{array}
\right.
\end{equation}

we are going to prove that the measure $\mu$ satisfy
\begin{equation}\label{massd}
\mu(J)\leq C|J|^d
\end{equation}
for any interval $J\subset [0,1]$, where $C$ is a positive constant
independent of $J$. We do this as following two steps.

{\bf Step 1.}  Suppose that $J$ is a component interval of $f(E_k)$,
For every $i, 0\leq i\leq k,$ let $J_i$ be the component interval of
$f(E_i)$ such that
\begin{equation}\label{sequence}
J=J_k\subset J_{k-1}\subset \cdot\cdot\cdot J_1\subset J_0=[0,1]
\end{equation}
By the definition of $\mu$, we have
$$
\frac{\mu(J)}{|J|^d}=\frac{1}{||J_{k-1}||_d}\frac{|J_{k-1}|^d}{||J_{k-2}||_d}\cdot\cdot\cdot
\frac{|J_{1}|^d}{||J_{0}||_d}=\frac{|J_{k-1}|^d}{||J_{k-1}||_d}\cdot\cdot\cdot
\frac{|J_{1}|^d}{||J_{1}||_d} \frac{|J_{0}|^d}{||J_{0}||_d}.
$$
Let
\begin{equation}
r_i=\frac{||J_{i}||_d}{|J_{i}|^d},\quad i=0,1,2,\cdot\cdot\cdot,k-1.
\end{equation}
So the above equality can be rewritten as
\begin{equation}
\frac{\mu(J)}{|J|^d}= (\prod_{i=1}^kr_{i-1})^{-1}.
\end{equation}
In order to prove (\ref{massd}), it suffices to show
\begin{equation}
\lim_{k\rightarrow \infty}\prod_{i=1}^kr_{i-1}=\infty.
\end{equation}

Given an $i$, $1\leq i\leq k$, we are going to estimate $r_{i-1}$. Let
$J_{i-1}$ be the component interval of $f(E_{i-1})$ in the sequence (\ref{sequence}).
Let $J_{i1},J_{i2},\cdot\cdot\cdot,J_{in_{i}}$ be the
 $n_i$ component intervals of $f(E_i)$ lying in $J_{i-1}$. Recall that $J_i\subset J_{i-1}$ is a component interval of
$f(E_i)$. So there must exist $1\leq i_0\leq n_i$ such that $J_i=J_{ii_0}$. Let
$G_{i0},G_{i1},\cdot\cdot\cdot,G_{in_i}$ be the $n_i+1$ gaps in the $J_{i-1}$. Put
$$
I_{i-1}=f^{-1}(J_{i-1}),\quad I_i=f^{-1}(J_i)=f^{-1}(J_{ii_0})\quad \mathrm{and}
\quad I_{ij}=f^{-1}(J_{ij}),
$$
for $j=1,2,\cdot\cdot\cdot,n_i.$
Then $I_{i1},\cdot\cdot\cdot,I_{in_i}$ are component intervals
of $E_i$ lying in the component interval $I_{i-1}$ of $E_{i-1}$.
Since $f$ is $M$-quasisymmetric, it follows Lemma
\ref{mainl2} and the construction of $E$ that
\begin{equation}\label{ineq14}
\frac{|G_{ij}|}{|J_{i-1}|}\leq 4(\frac{|f^{-1}(G_{ij})|}{|f^{-1}(J_{i-1})|})^p\leq 4e_i^p,\quad
j=0,1,2,\cdot\cdot\cdot,n_i,
\end{equation}
where $e_i=\max_{0\leq l\leq n_i}{e_{i,l}}$
and that
\begin{equation}\label{ineq15}
\frac{|J_{ij}|}{|J_{i-1}|}\geq
(1+M)^{-2}(\frac{|I_{ij}|}{|I_{i-1}|})^q=(1+M)^{-2}c_i^q.
\end{equation}
Here $p,q$ are numbers defined in Lemma \ref{mainl2}. The
inequality (\ref{ineq14}) yields
\begin{equation}\label{E14}
\frac{|J_{i1}|+\cdot\cdot\cdot+|J_{in_i}|}{|J_{i-1}|}=
\frac{|J_{i-1}|-|G_{i0}|-\cdot\cdot\cdot-|G_{in_i}|}{|J_{i-1}|}\geq
1-4(n_i+1)e_i^p.
\end{equation}
From inequality (\ref{ineq15}), we have
\begin{equation}\label{E20}
\begin{split}
r_{i-1}&\;=\frac{|J_{i1}|^d+\cdot\cdot\cdot+|J_{in_i}|^d}{|J_{i-1}|^d}\\
&\;\geq n_i(\frac{|J_{ij}|}{|J_{i-1}|})^{d}\\
&\;\geq \frac{n_i}{(1+M)^{2d}}c_i^{dq}.
\end{split}
\end{equation}
Let
$$
S(k,p)=\{i:1\leq i\leq k,e^p_i\leq \min\{a,|I_i|^p\}
$$
 where $a=1-\sqrt[4N+4]{\frac{4N+4}{4N+5}}$, where $N=1+\sup_{l} n_l$. 
 Since $\eta_{i,l}\leq e_{i,l}$. Thus, If $i\in S(k,p)$ we have
\begin{equation}
\begin{split}
c_i=\frac{|I_{ij}|}{|I_{i-1}|}&\;=\frac{|I_{ij}|}{n_i|I_{ij}|+\sum_{l=0}^{n_i}\eta_{i,l}}\\
&\;\geq\frac{|I_{ij}|}{n_i|I_{ij}|+(n_i+1)\eta_i}\\
&\;\geq\frac{1}{2n_i+1}\\
&\;\geq\frac{1}{2N}
\end{split}
\end{equation}
for $j=1,\cdot\cdot\cdot,n_i$, where $\eta_i=\max_{0\leq l\leq n_i}{\eta_{i,l}}$.

From the conclusion (3) of Lemma \ref{Mlema}, we obtain
 \begin{equation}\label{eq0000}
 \lim_{k\rightarrow \infty}\frac{\sharp S(k,p)}{k}=1.
 \end{equation}


Then follows from the left hand inequality of (\ref{numpq1}) that
$$
1\geq\frac{|J_{ij}|}{|J_i|}
=\frac{|f(I_{ij})|}{|f(I_i)|}\geq(1+M)^{-2}(\frac{|I_{ij}|}{|I_{i-1}|})^{q}\geq
A
$$
for $j=1,2,\cdot\cdot\cdot,n_i$, where
$A=\frac{(1+M)^{-2}}{(2N)^{q}}$. Therefore,
\begin{equation}\label{E16}
\begin{split}
\frac{|J_i|^d+|J_{i1}|^d+\cdot\cdot\cdot+|J_{in_i}|^d}
{(|J_i|+|J_{i1}|+\cdot\cdot\cdot+|J_{in_i}|)^d}
&\;=\frac{1+x_1^d+\cdot\cdot\cdot+x_{n_i}^d}{(1+x_1+\cdot\cdot\cdot+x_{n_i})^d}\\
 &\;\geq (1+A)^{1-d},
\end{split}
\end{equation}
where $x_j=\frac{|J_{ij}|}{|J_i|}\in [A,1]$.

Note that the equality (\ref{E14}) and (\ref{E16}), for any $i\in S(k,p)$ we obtain
\begin{equation}
\begin{split}
r_{i-1}&\;=\frac{|J_i|^d+|J_{i1}|^d+\cdot\cdot\cdot+|J_{in_i}|^d}{|J_{i-1}|^d}\\
&\;=\frac{|J_i|^d+|J_{i1}|^d+\cdot\cdot\cdot+|J_{in_i}|^d}
{(|J_i|+|J_{i1}|+\cdot\cdot\cdot+|J_{in_i}|)^d}\frac{(|J_i|+|J_{i1}|+\cdot\cdot\cdot+|J_{in_i}|)^d}
{|J_{i-1}|^d}\\
&\;\geq \alpha_2(1-4(n_i+1)e_i^p)^d,
\end{split}
\end{equation}
where $\alpha_2=(1+A)^{1-d}>1$.

Since
 $$
 1-mx\geq (1-x)^{m+1}
 $$
for all $x\in (0,1-\sqrt[m]{\frac{m}{m+1}})$, so we have
$$
1-4mx\geq (1-x)^{4m+1}
$$
for all $x\in (0,a)$where $a=1-\sqrt[4N+4]{\frac{4N+4}{4N+5}}$ and all positive inters $m\leq N$.

Note that $n_i<N$ and $e_i^p\in (0,a)$ for all $i\in S(k,p)$, thus we obtain
\begin{equation}\label{E21}
r_{i-1}\geq \alpha_2(1-e_i^p)^{(4n_i+4)d}
\end{equation}
Using the estimate (\ref{E20}) and (\ref{E21}), we obtain
\begin{equation}
\begin{split}
\prod_{i=1}^kr_{i-1}&\;\geq \prod_{i \not\in
S(k,p)}\frac{n_ic_i^{dq}}{(1+M)^{2d}}
\prod_{i\in S(k,p)}\alpha_2(1-4(n_i+1)e_i^p)^d\\
&\;\geq\prod_{i \not\in
S(k,p)}\frac{n_ic_i^{dq}}{(1+M)^{2d}}\prod_{i\in
S(k,p)}\alpha_2(1-e_i^p)^{(4n_i+4)d}\\
&\;=\alpha_2^{\sharp S(k,p)}[(1+M)^{-2d}]^{k-\sharp S(k,p)}\prod_{i \not\in
S(k,p)}n_ic_i^{dq}\prod_{i\in
S(k,p)}(1-e_i^p)^{(4n_i+4)d}.
\end{split}
\end{equation}


If $q=1$, since $n_ic_i\leq 1$ then we have
$$
\prod_{i \not\in S(k,p)}n_ic_i^{dq}=\prod_{i \not\in S(k,p)}n_ic_i^{d}\geq\prod_{i \not\in S(k,p)}n_ic_i\geq\prod^k_{i =1}n_ic_i=N_k\delta_k.
$$
If $q>1$, we have
\begin{equation}
\begin{split}
\prod_{i \not\in
S(k,p)}n_ic_i^{dq}&\;=\prod^k_{i\not\in S(k,p)
}(n_ic_i)^{dq}n_i^{1-dq}\geq\prod^k_{i =1}(n_ic_i)^{dq}\prod_{i \not\in
S(k,p)}n_{i}^{1-dq}\\
&\;=\prod^k_{i =1}(n_ic_i)^{dq}\prod_{i \not\in
S(k,p)}n_{i}^{1-dq}
\geq (N_k\delta_k)^{dq}\prod_{i \not\in
S(k,p)}N^{1-dq}\\
&\;=(N_k\delta_k)^{dq}(N^{1-dq})^{k-{\sharp S(k,p)}}
\end{split}
\end{equation}
for $d\in (1/q,1)$.

Let
\begin{equation}\label{eq1}
\xi_k=\alpha_2^{\sharp S(k,p)}((1+M)^{-2d})^{k-\sharp S(k,p)}(N_k\delta_k)^{dq}(N^{1-dq})^{k-\sharp S(k,p)}
\end{equation}
and
$$
 \zeta_k=\prod_{i\in
S(k,p)}(1-e_i^p)^{(4n_i+4)d}.
$$
Thus, we have
\begin{equation}\label{eq0}
\prod_{i=1}^kr_{i-1}\geq  \xi_k\zeta_k.
\end{equation}
It is obvious that
\begin{equation}\label{eq2}
\lim_{k\rightarrow \infty}\xi_k^{1/k}=\alpha_2>1.
\end{equation}
due to the conclusion (1) of Lemma \ref{Mlema} and the equality (\ref{eq0000}) .
On the other hand, since $\log(1-x)\geq -2x$ when $0<x<1$, the conclusion (2) of Lemma \ref{Mlema}, we obtain
\begin{equation}
\begin{split}
\frac{1}{k}\log\zeta_k&\;=\frac{1}{k}\log\prod_{i\in
S(k,p)}(1-e_i^p)^{(4n_i+4)d}\\
&\;=\frac{1}{k}\sum_{i\in
S(k,p)}\log(1-e_i^p)^{(4n_i+4)d}\\
&\;=\frac{1}{k}\sum_{i\in
S(k,p)}(4n_i+4)d\log(1-e_i^p)\\
&\;\geq \frac{(4N+4)d}{k}\sum_{i\in
S(k,p)}\log(1-e_i^p)\\
&\;\geq -2\frac{(4N+4)d}{k}\sum_{i\in
S(k,p)}e_i^p\\
&\;\geq -2\frac{(4N+4)d}{k}\sum_{i=1}^ke_i^p\rightarrow 0.
\end{split}
\end{equation}
as $k\rightarrow\infty$. This show that
\begin{equation}\label{eq3}
\lim_{k\rightarrow \infty}\zeta_k^{1/k}=1.
\end{equation}
From $(\ref{eq0}),(\ref{eq2}),(\ref{eq3})$, we obtain
$$
\liminf_{k\rightarrow \infty}(\prod_{i=1}^kr_{i-1})^{1/k}\geq \alpha_2>1.
$$
This implies
$$
\lim_{k\rightarrow \infty}(\prod_{i=1}^kr_{i-1})=\infty.
$$

{\bf Step 2.} Let $J\subset [0,1]$ be any interval. For such $J$, let $k$ be the unique positive inter such that
$$
\delta_k\leq |f^{-1}(J)|\leq \delta_{k-1},
$$
where $\delta_k$ denotes the lengthen of component intervals of $E_k$. Then the set $f^{-1}(J)$ meets at most two component intervals of $E_{k-1}$ and hence at most $2n_{k+1}$ component intervals of $E_k$. Thus, the set $J$ meets at most $2n_{k+1}$ component intervals of $f(E_k)$.

Let $J_1,J_2,\cdot\cdot\cdot, J_l,l\leq 2n_{k+1}$, be those component intervals of $f(E_k)$ meeting $J$. Using the conclusion of step 1. we obtain
\begin{equation}
\mu(J)\leq \sum_{i=1}^l \mu(J_i)\leq C\sum_{i=1}^l |J_i|^d.
\end{equation}
 Since $\delta_k\leq |f^{-1}(J)|$, we obtain
 $$
 f^{-1}(J_i)\subset 3f^{-1}(J), \quad i=1,2,3\cdot\cdot\cdot l,
 $$
 where $3f^{-1}(J)$ denote the interval of lengthen $3|f^{-1}(J)|$ concentric with $f^{-1}(J)$. Thus we obtain
 $$
 |J_i|\leq f(3f^{-1}(J))\leq K|J|, \quad i=1,2,3\cdot\cdot\cdot l,
 $$
where $K$ is a positive constant depending on $M$ only. This together with gives
$$
\mu(J)\leq ClK^d|J|^d\leq 2NCK^d|J|^d.
$$
This show that (\ref{massd}).

By Lemma (\ref{MSSLemma}), it follows from  that $\dim_Hf(E)\geq d$ for  $d$. As $d$ could be chosen as closed to 1 as one would. Hence $\dim_Hf(E)=1$. This completes the proof of Theorem \ref{MTHM}.\qed

{\bf Acknowledgments}. I would like to thank my advisor Professor Qiu Weiyuan for introducing me to the theory of fractal geometry.

\thebibliography{10}
\bibitem{KF}K.J.Faconer, Fractal Geometry: Mathematical Foundations and Applications, John Wile Sons (1990).
%

\bibitem{HH}H.Hakobyan,  Cantor sets minimal for quasisymmetric maps, J.Contemp. Math.Anal. {\bf 41(2)}, 2006, 5-13.

\bibitem{MS}M.D.Hu, S.Y.Wen, Quasisymmetrically minimal uniform cantor sets, Topology and its
Applications, {\bf 155}, 2008, no.6,515-521.

\bibitem{JH}J.Heinonen, Lectures on Analysis in Metric Spaces. Universitext. Springer-Verlag, New York, 2001.

\bibitem{PP}P.Pansu,  Dimension confrome et sph\`{e}re \`{a} l'infini des vari\'{e}t\'{e}s \`{a} courbure n\'{e}gative, Ann.Acad.Sci.Fenn. {\bf14(2)}14 177-212.

\bibitem{ZW}Z.Y.Wen, J.Wu, Hausdorff dimension of homogeneous perfect sets, Acta Math.Hungar.,{\bf107},2005,35-44.

\bibitem{ZYW}Z.Y.Wen, Mathematical Foundations of Fractal geometry, Shanghai Scientific and Technological Education Publishing House, 2000.

\bibitem{WW}X.Y.Wang, J.Wu, Packing dimensions of homogeneous perfect
    sets, Acta Math.Hungar., {\bf118(1-2)},2008,29-39.

\bibitem{JMW}J.M.Wu, Null sets for doubling and dyadic doubling measures, Ann.Acad.Sci.Fenn.Math.{\bf18}, 1993,77-91.
\end{document}